\documentclass[twoside,twocolumn]{article}

\usepackage{blindtext} 

\usepackage[sc]{mathpazo} 
\usepackage[T1]{fontenc} 
\usepackage{microtype} 

\usepackage[english]{babel} 

\usepackage[hmarginratio=1:1,top=28mm,columnsep=20pt,bottom=64pt,left=80pt]{geometry} 
\usepackage[hang, small,labelfont=bf,up,textfont=it,up]{caption} 
\usepackage{booktabs} 

\usepackage{lettrine} 

\usepackage{enumitem} 
\setlist[itemize]{noitemsep} 

\usepackage{abstract} 

\usepackage{titlesec} 
\renewcommand\thesection{\arabic{section}} 
\renewcommand\thesubsection{\arabic{section}.\arabic{subsection}} 
\titleformat{\section}[block]{\Large\scshape\centering}{\Large\thesection}{1em}{} 
\titleformat{\subsection}[block]{\large}{\thesubsection}{1em}{} 

\usepackage{fancyhdr} 
\usepackage{titling} 

\usepackage{hyperref} 

\usepackage{amsmath}
\usepackage{amsfonts}
\usepackage{graphicx}
\usepackage{float}

\usepackage{xcolor}
\hypersetup{
	colorlinks,
	linkcolor={black},
	citecolor={black},
	urlcolor={black}
}

\pretitle{\begin{center}\LARGE\bfseries} 
\posttitle{\end{center}} 
\title{The Inscribed Angle Theorem for the Hyperbola} 
\author{%
\textsc{Jack Williams}
\\[1ex] 
\normalsize University of Cambridge \\ 
\normalsize \href{mailto:jack.williams@damtp.cam.ac.uk}{jack.williams@damtp.cam.ac.uk} 
}
\date{\today} 
\renewcommand{\maketitlehookd}{%
\begin{abstract}
\noindent The inscribed angle theorem, a famous result about the angle subtended by a chord within a circle, is well known and commonly taught in school curricula.  In this paper, we present a generalisation of this result (and other related circle theorems) to the rectangular hyperbola.  The notion of angle is replaced by pseudo-angle, defined via the Minkowski inner product.  Indeed, in Minkowski space, the unit hyperbola is the set of points a unit metric distance from the origin, analogous to the Euclidean unit circle.  While this is a result of pure geometrical interest, the connection to Minkowski space allows an interpretation in terms of special relativity where, in the limit $c\to\infty$, it leads to a familiar result from non-relativistic dynamics.  This non-relativistic result can be interpreted as an inscribed angle theorem for the parabola, which we show can also be obtained from the Euclidean inscribed angle theorem by taking the limit of a family of ellipses ananlogous to the non-relativistic limit $c\to\infty$. 
This simple result could be used as a pedagogical example to consolidate understanding of pseudo-angles in non-Euclidean spaces or to demonstrate the power of analytic continuation.
\end{abstract}
}

\begin{document}

\maketitle


\section{Introduction}

By replacing $t$ with $it$, many results and techniques familiar in the Euclidean setting can be applied in physically meaningful Lorentzian settings \cite{Wick:1954eu}.  Making use of these Wick rotations, we obtain a generalisation of the inscribed angle theorem to rectangular hyperbolae in Minkowski space.  The hyperbola itself is analogous to a Euclidean circle and is related to it by a Wick rotation of the plane coordinates.  The Euclidean angle is replaced by the Minkowski pseudo-angle and the proof of the hyperbolic inscribed angle theorem is similar to a proof of the Euclidean result, but with trigonometric functions replaced by hyperbolic trigonometric functions.

While this is a result of pure geometrical interest, the connection to Minkowski space allows an interpretation in terms of special relativity.  Observing that hyperbolae in Minkowski space occur naturally as the trajectory of observers moving with constant proper acceleration, we provide a physical interpretation of the inscribed angle theorem in terms of an observer who ejects two particles while moving along this trajectory.  The observer subsequently collides with both particles. The inscribed angle theorem implies that the proper time between these two events depends only on the relative velocity of the two particles and not on the relative velocity between the particles and the observer.  This is the relativistic analogue of a similar fact familiar from non-relativistic kinematics and, of course, the non-relativistic result can be obtained from the relativistic one by taking the limit $c\to\infty$.

Since the trajectory of a non-relativistic observer moving with constant acceleration is a parabola, this limit leads to an inscribed angle theorem for the parabola.  Since the parabola is the singular conic between families of hyperbolae and ellipses, it is natural that the inscribed angle theorem should arise by taking this limit of a family of hyperbolae.  We show that the same result can also be obtained from the Euclidean inscribed angle theorem by taking a limit of a family of ellipses.

This result is accessible to early undergraduates and could be presented as a pedagogical example to build intuition around the analogy between Euclidean angles and rapidity in special relativity, or Euclidean rotations and Lorentz transformations.  It also demonstrates the power of analytic continuation.

\section{Inscribed Angle Theorem on the Unit Hyperbola}
\label{inscribedangle}

\subsection{Inscribed angle theorem in 2D Euclidean space}
In 2D Euclidean space, the inscribed angle theorem states that given any three points $P_0,P_1,P_2$ on the unit circle, the lines $P_0P_1$ and $P_0P_2$ intersect at an acute angle independent of the position of $P_0$ \cite{euclid}.  Further, if the points are represented as $P_i=(\cos(\theta_i),\sin(\theta_i))$, where $0\leq\theta_i<2\pi$, then the acute angle between $P_0P_1$ and $P_0P_2$ is $\frac{1}{2}|\theta_1-\theta_2|$ (or $\pi-\frac{1}{2}|\theta_1-\theta_2|$, whichever is acute).

This follows because the gradients of the chords $P_0P_i$ are 
\begin{align}
\frac{\sin\theta_0-\sin\theta_i}{\cos\theta_0-\cos\theta_i}&=\resizebox{.25\textwidth}{!}{$\frac{2\cos\left(\frac{\theta_0+\theta_i}{2}\right)\sin\left(\frac{\theta_0-\theta_i}{2}\right)}{2\sin\left(\frac{\theta_0+\theta_i}{2}\right)\sin\left(\frac{\theta_0-\theta_i}{2}\right)}$}\nonumber
\\&=\cot\left(\frac{\theta_0+\theta_i}{2}\right).
\end{align}

Given that the angle between them satisfies
\begin{align}
&\cos^2(\theta)=\frac{(x\cdot y)^2}{x^2y^2},
\end{align}

where $x$ and $y$ are tangent to these chords, we find  
\begin{align}
\cos^2(\theta)&=
 \resizebox{.33\textwidth}{!}{$\frac{\left(1+\cot\left(\frac{\theta_0+\theta_1}{2}\right)\cot\left(\frac{\theta_0+\theta_2}{2}\right)\right)^2}{\left(1+\cot^2\left(\frac{\theta_0+\theta_1}{2}\right)\right)\left(1+\cot^2\left(\frac{\theta_0+\theta_1}{2}\right)\right)}$}\nonumber\\&=\cos^2\left(\frac{\theta_1-\theta_2}{2}\right).
\end{align}

It follows immediately that $\theta$ is independent of $\theta_0$ and depends only the the difference between $\theta_1$ and $\theta_2$.

\subsection{Pseudo-angles in 2D Minkowski space}

It can be checked directly that the inner product on $\mathbb{R}^{1,1}$ given by $x\cdot y=-x_0y_0+x_1y_1$ satisfies a form of the Cauchy-Schwarz inequality: for any non-null vectors $x$ and $y$, $(x\cdot y)^2\geq x^2y^2$.  Using this, we can define the pseudo-angle between two vectors in $\mathbb{R}^{1,1}$ via $\cosh^2(\theta)=\frac{(x\cdot y)^2}{x^2y^2}$ whenever $x$ and $y$ are either both spacelike or both timelike \cite{nabel}.  A similar construction exists for the pseudo-angle between a spacelike and a timelike vector.

In the case that $x$ and $y$ are both timelike, this pseudo-angle is, in fact, the relative rapidity of observers $A$ and $B$ moving with velocities $x$ and $y$.  Indeed, if $B$ moves with velocity vector parallel to $(1,v)$ in the rest frame of $A$, then 
\begin{align}
\cosh^2(\theta)&=\frac{(x\cdot y)^2}{x^2y^2}=\frac{1}{1-v^2} \nonumber
\\&\implies v=\pm\tanh(\theta).
\end{align}

\subsection{Inscribed angle theorem in 2D Minkowski space}
In 2D Minkowski space, the analogue of the unit circle is the hyperbola $-x_0^2+x_1^2=1$, consisting of points a unit Minkowski distance from the origin.  Let $P_0,P_1,P_2$ be any three points on the branch of this hyperbola with $x_1>0$ with coordinates $P_i=(\sinh(\theta_i),\cosh(\theta_i))$.  Observe that the chords $P_0P_i$ must be timelike.

\begin{figure}[H]
	\centering
	\includegraphics[scale=.45]{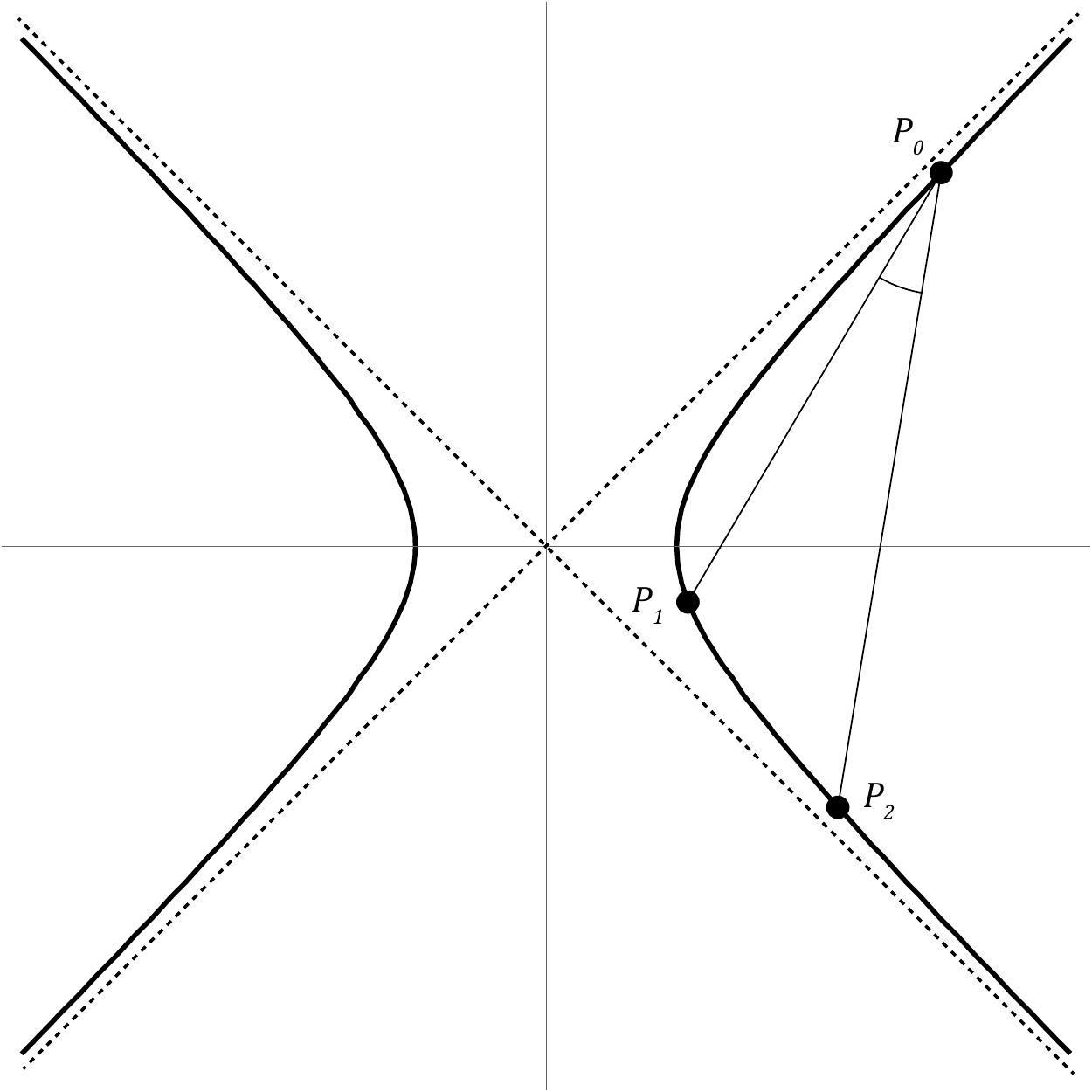}
	\caption{Inscribed angle theorem in Minkowski space}
\end{figure}

In a similar way to the Euclidean case, we find that the gradient of the line $P_0P_i$ is $\coth\left(\frac{\theta_0+\theta_i}{2}\right)$ and hence that the pseudo-angle between these chords is given by

\begin{align}
\cosh^2(\theta)&= \resizebox{.35\textwidth}{!}{$\frac{\left(\coth\left(\frac{\theta_0+\theta_1}{2}\right)\coth\left(\frac{\theta_0+\theta_2}{2}\right)-1\right)^2}{\left(\coth^2\left(\frac{\theta_0+\theta_1}{2}\right)-1\right)\left(\coth^2\left(\frac{\theta_0+\theta_1}{2}\right)-1\right)}$}\nonumber\\&=\cosh^2\left(\frac{\theta_1-\theta_2}{2}\right).
\end{align}

It follows again that $\theta$ does not depend on the position of $P_0$.

\subsection{Analogues of other circle theorems in Minkowski space}
There are similar analogues of other familiar Euclidean results.  For example, on the circle, the angle $P_1P_0P_2$ is half the angle $P_1OP_2$ \cite{euclid}.  This holds true in Minkowski space since the pseudo-angle between the spacelike chords $OP_1$ and $OP_2$ is given by 
\begin{align}
&\cosh^2(\theta)=\frac{(x\cdot y)^2}{x^2y^2}\nonumber\\&=\left(\cosh\left(\theta_1\right)\cosh\left(\theta_2\right)-\sinh\left(\theta_1\right)\sinh\left(\theta_2\right)\right)^2\nonumber\\&=\cosh^2\left(\theta_1-\theta_2\right).
\end{align}
and is hence twice $\frac{\theta_1-\theta_2}{2}$, the pseudo-angle between $P_0P_1$ and $P_0P_2$.

\begin{figure}[H]
	\centering
	\includegraphics[scale=.45]{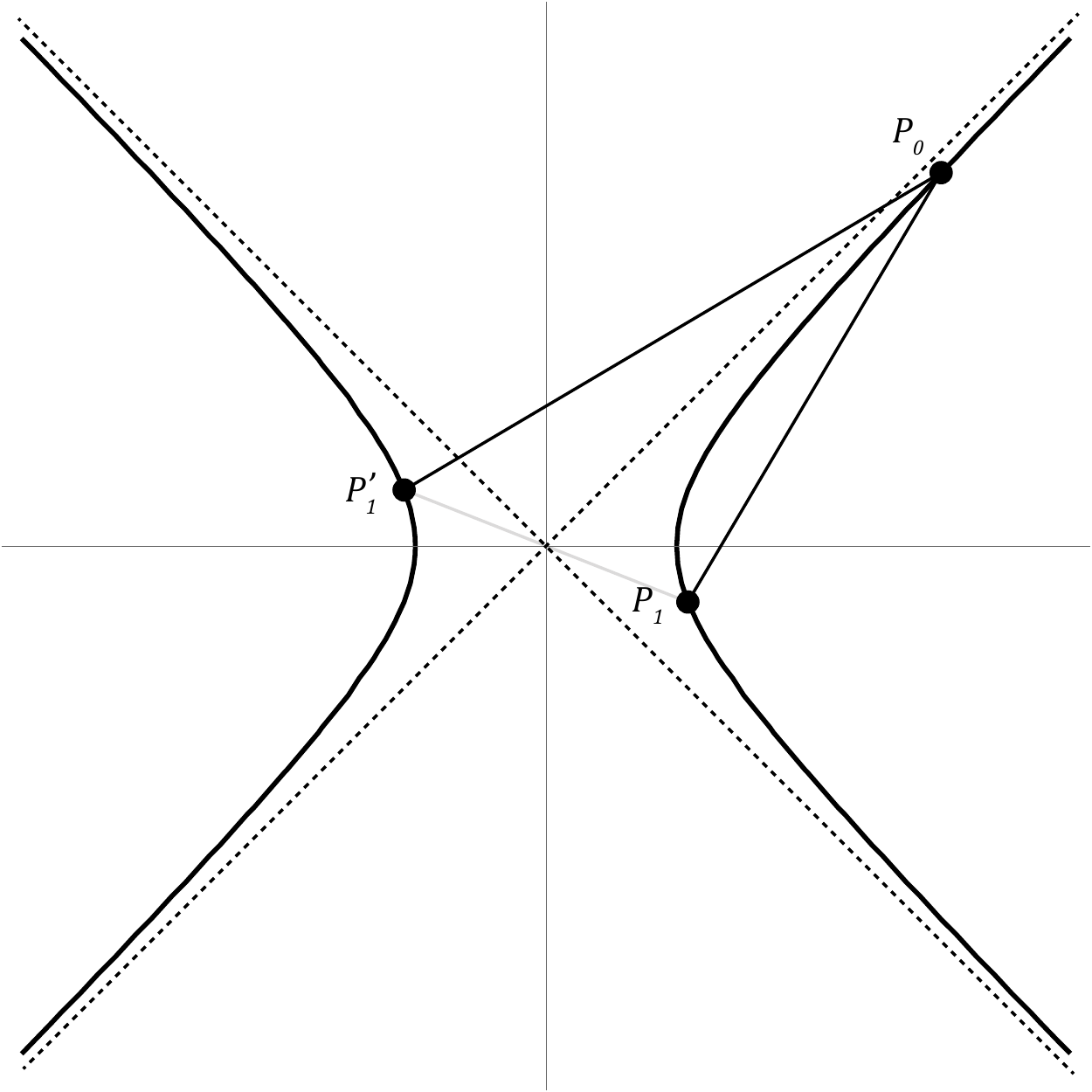}
	\caption{Thales' theorem in Minkowski space}
\end{figure}

Furthermore, the analogue of Thales' theorem, which states that the angle subtended by a diameter is a right angle \cite{euclid}, is the result that if $P_1'=-P_1$ is the point in Minkowski space obtained from $P_1$ by a reflection in the origin, then the chords $P_0P_1$ and $P_0P_1'$ are orthogonal.  This holds because the gradients of $P_0P_1$ and $P_0P_1'$ are $\coth\left(\frac{\theta_0+\theta_1}{2}\right)$ and $\tanh\left(\frac{\theta_0+\theta_1}{2}\right)$ respectively.  Their tangent vectors are orthogonal since
\begin{align}
x\cdot y=\resizebox{.25\textwidth}{!}{$\coth\left(\frac{\theta_0+\theta_1}{2}\right)\tanh\left(\frac{\theta_0+\theta_1}{2}\right)$}-1=0.
\end{align}

\subsection{Relationship between the inscribed angle theorem and isometries}
In the traditional statement of the inscribed angle theorem in Euclidean space, the points $P_1$ and $P_2$ are fixed and the point $P_0$ varies on the circle.  However, to obtain an alternative viewpoint one keeps $P_0$ fixed and rotates the two lines $P_0P_1$ and $P_0P_2$ together, keeping their relative angle fixed.  After the rotation, the new points of intersection, $P_1'$ and $P_2'$, of these lines with the circle are separated by the same relative angle as the original points $P_1$ and $P_2$.  In fact, the rotation taking $P_1$ and $P_2$ to $P_1'$ and $P_2'$ is an isometry of the Euclidean plane.

\begin{figure}[H]
	\centering
	\includegraphics[scale=.45]{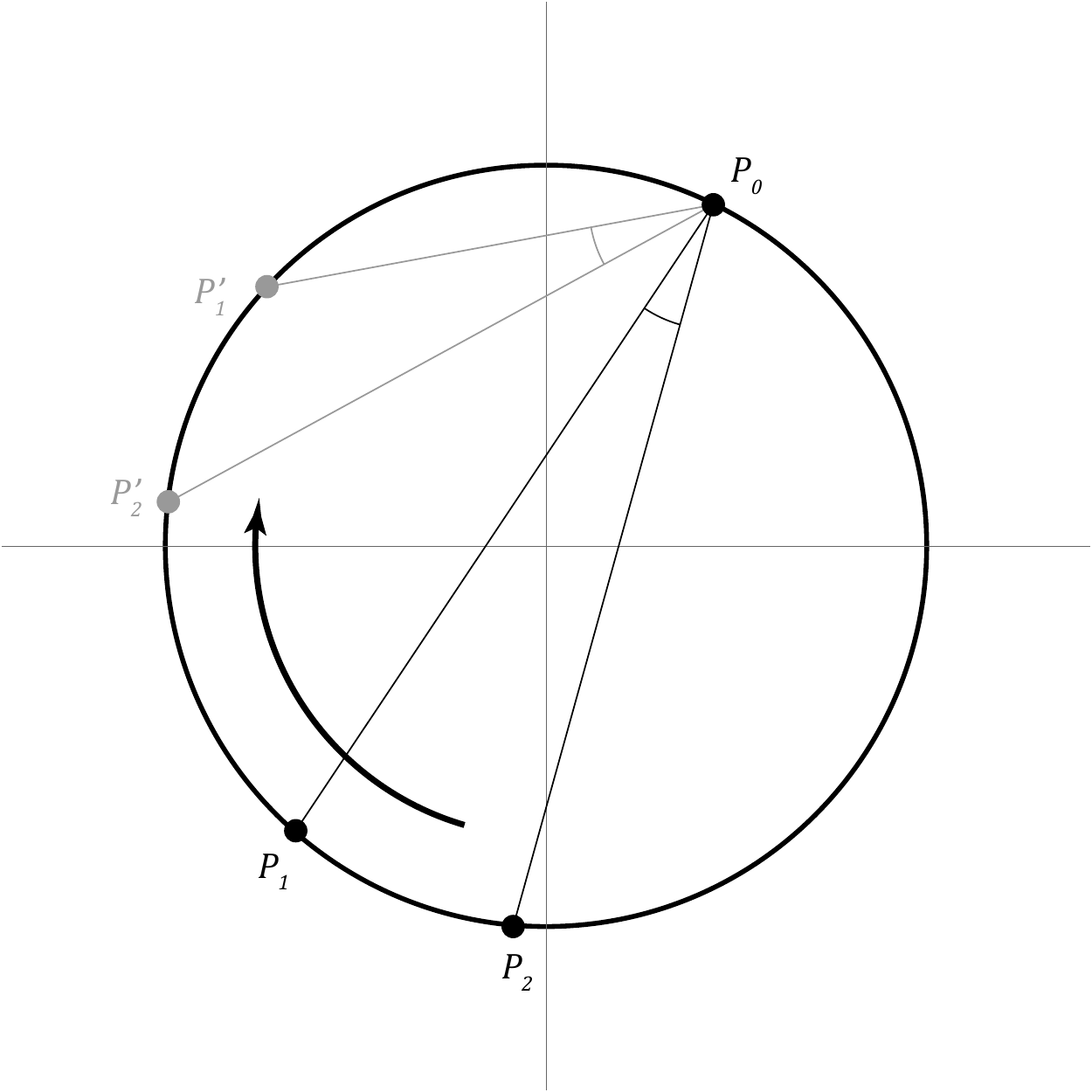}
	\caption{Alternative interpretation of Euclidean inscribed angle theorem}
\end{figure}

Likewise, in Minkowski space one can rotate the lines $P_0P_1$ and $P_0P_2$ together about $P_0$, keeping their relative pseudo-angle fixed. The pseudo-angular distance between the new points of intersection, $P_1'$ and $P_2'$, of these lines with the hyperbola remains $\frac{\theta_1-\theta_2}{2}$.  Since isometries of Minkowski space of the form 
\begin{align}
\left(\begin{matrix}x_0\\x_1\end{matrix}\right)\mapsto\left(\begin{matrix}\cosh\phi&\sinh\phi\\\sinh\phi&\cosh\phi\end{matrix}\right)\left(\begin{matrix}x_0\\x_1\end{matrix}\right)
\end{align}
map points $(\sinh(\theta_i),\cosh(\theta_i))$ to $(\sinh(\theta_i+\phi),\cosh(\theta_i+\phi))$, there is an isometry of Minkowski space mapping $P_1'$ and $P_2'$ back to $P_1$ and $P_2$.  This isometry relates the two perspectives on the inscribed angle theorem.

\begin{figure}[H]
	\centering
	\includegraphics[scale=.45]{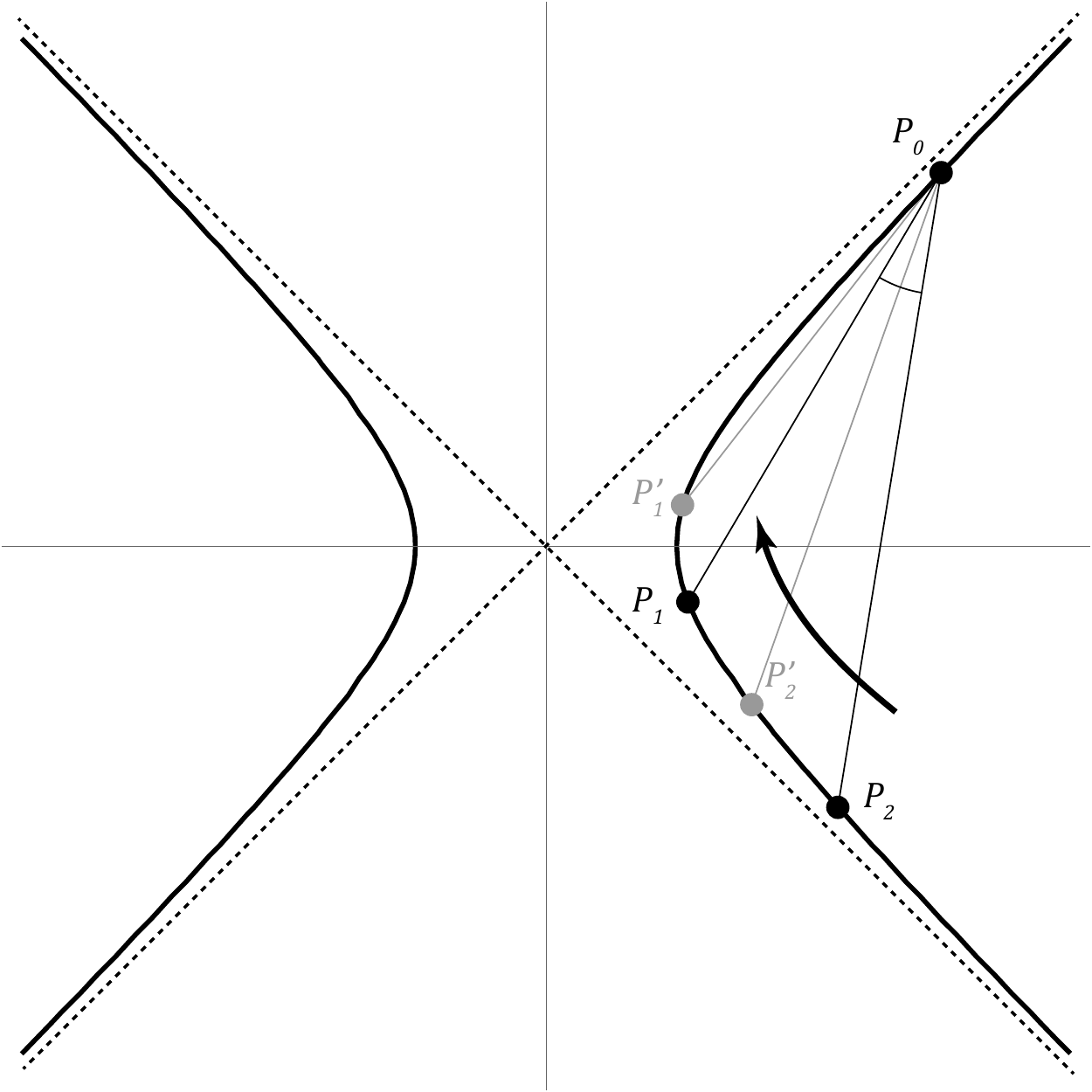}
	\caption{Alternative interpretation of hyperbolic inscribed angle theorem}
\end{figure}

\section{Relativistic interpretation}

\subsection{Constant proper acceleration}
An observer moving with constant unit proper acceleration in Minkowski space follows a rectangular hyperbola of the form $-x_0^2+x_1^2=1$ \cite{misner}.  Suppose this observer ejects two particles, $A$ and $B$ with fixed relative rapidity (i.e. separated by a fixed pseudo-angle), in the same direction as the observer's acceleration but with greater velocity.  Then at some future times, the observer will collide again with each of these objects.  
\begin{figure}[H]
	\centering
	\includegraphics[scale=.45]{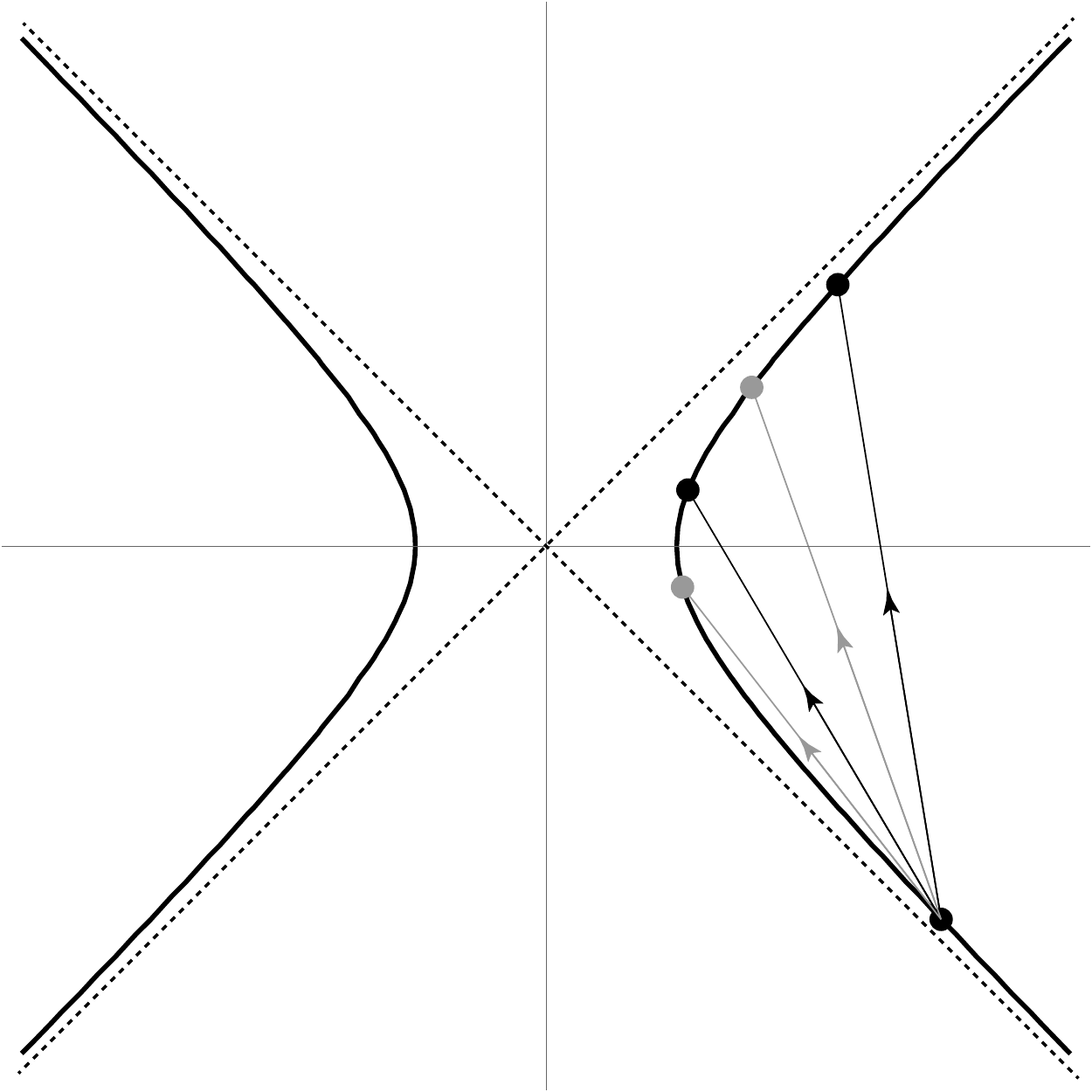}
	\caption{Kinematic interpretation of hyperbolic inscribed angle theorem}
\end{figure}

Since the position of the observer at proper time $\tau$ is $(\sinh(k(\tau-\tau_0)),\cosh(k(\tau-\tau_0)))$, the inscribed angle theorem (viewed from the perspective in which $P_0$ is fixed and the two lines rotate) implies that the proper time separating these two collisions is independent of the velocity with which $A$ was ejected by the observer.  It depends only on the relative velocity between $A$ and $B$.

\subsection{Connection with non-relativistic kinematics}
This is an analogue of another familiar fact from non-relativistic kinematics.  An observer with constant acceleration follows a trajectory of the form $x=\frac{1}{2}at^2+ut+x_0$.  If, at time 0, two particles are ejected with velocities $\alpha$ and $\beta$, then these particles follow trajectories $x=\alpha t+x_0$ and $x=\beta t+x_0$.  The difference between the times at which the observer collides with these particles is $\frac{2}{a}(\alpha-\beta)$, which depends only on the relative velocity of the particles.  

A corresponding inscribed angle theorem for the parabola arises from this non-relativistic result, where the pseudo-angle between two lines $P_0P_1$ and $P_0P_2$ is proportional to the difference in the $t$-coordinates of $P_1$ and $P_2$ and is independent of $P_0$.  The inscribed angle theorem states that this angle measure is constant as $P_0$ varies on the parabola.

\begin{figure}[H]
	\centering
	\includegraphics[scale=.45]{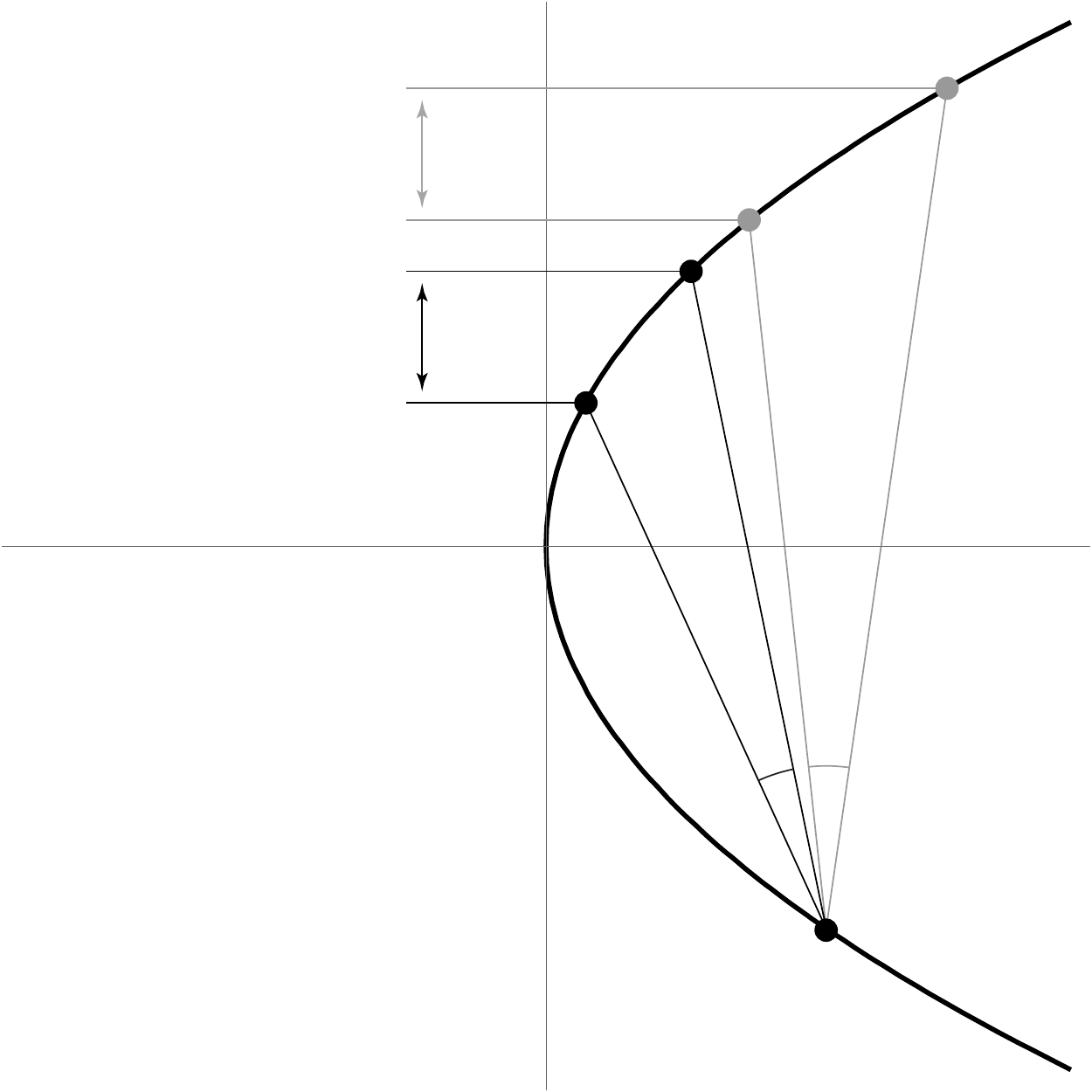}
	\caption{Inscribed angle theorem for parabola}
\end{figure}

Naturally, this non-relativistic inscribed angle theorem for the parabola can be obtained from the inscribed angle theorem for the hyperbola by taking the limit $c\to\infty$.  In order to do this, it is helpful to reintroduce dimensionful quantities.  The trajectory of a relativistic observer moving with constant proper acceleration $a$ in the positive $x$-direction and passing through the origin is given by the hyperbola 
\begin{align}
\resizebox{.42\textwidth}{!}{$x=\frac{c^2}{a} \left(\cosh\left( \frac{a\theta}{c} \right)-1\right),\ t=\frac{c}{a} \sinh\left( \frac{a\theta}{c} \right)$}
\end{align}
where $\theta$ is the observer's proper time.  In the non-relativistic limit, these become 
\begin{align}
\label{parabola}
x=\frac{1}{2}a\theta^2,\ t= \theta
\end{align}
and we recover the non-relativistic parabola.  The gradient of a chord $P_0P_i$ is $\frac{2}{a(\theta_0+\theta_i)}$ where $P_i=\left(\theta_i,\frac{1}{2}a\theta_i^2\right)$.

Using the metric $ds^2=-c^2dt^2+dx^2$ to compute inner products, the angle between the chords $P_0P_1$ and $P_0P_2$ is given by 
\begin{align}
\cosh^2\left(\theta\right)=\resizebox{.29\textwidth}{!}{$\frac{\left( \frac{4}{a^2(\theta_0+\theta_1)(\theta_0+\theta_2)}-c^2 \right)^2}{\left( \frac{4}{a^2(\theta_0+\theta_1)^2}-c^2 \right)\left( \frac{4}{a^2(\theta_0+\theta_2)^2}-c^2 \right)}$}
\end{align}
and so $\cosh^2\left(\theta\right)\to1$ as $c\to\infty$.  Using $\cosh^2\left(\theta\right)\sim1+\theta^2$, we find $\theta^2\sim\frac{a^2}{4c^2}(\theta_1-\theta_2)^2$, recovering the result that the angle between $P_0P_1$ and $P_0P_2$ is proportional to the difference in the $t$-coordinates of $P_1$ and $P_2$ and is independent of $P_0$.

\subsection{The parabolic limit of the Euclidean inscribed angle theorem}

The inscribed angle theorem for the parabola can also be obtained in a suitable limit from the Euclidean inscribed angle theorem.  Consider the ellipse $\left(x-\frac{c^2}{a}\right)^2+c^2t^2=\frac{c^4}{a^2}$ passing through the origin.  In the limit $c\to\infty$, this becomes the parabola $x=\frac{1}{2}at^2$.  It can be parametrised as 
\begin{align}
\resizebox{.41\textwidth}{!}{$x=\frac{c^2}{a} \left(1-\cos\left( \frac{a\theta}{c} \right)\right),\ t=\frac{c}{a} \sin\left( \frac{a\theta}{c} \right)$}
\end{align}
which is the limit $c\to\infty$ takes the form (\ref{parabola}).  The gradient of the chord $P_0P_i$ is $\frac{2}{a(\theta_0+\theta_i)}$ where $P_i=\left(\theta_i,\frac{1}{2}a\theta_i^2\right)$.

Using the metric $ds^2=dx^2+c^2dt^2$, in which the ellipse is a metric circle, the angle between $P_0P_1$ and $P_0P_2$ is given by 
\begin{align}
\cos^2\left(\theta\right)=\resizebox{.29\textwidth}{!}{$\frac{\left( \frac{4}{a^2(\theta_0+\theta_1)(\theta_0+\theta_2)}+c^2 \right)^2}{\left( \frac{4}{a^2(\theta_0+\theta_1)^2}+c^2 \right)\left( \frac{4}{a^2(\theta_0+\theta_2)^2}+c^2 \right)}$}
\end{align}
which implies $\theta^2\sim\frac{a^2}{4c^2}(\theta_1-\theta_2)^2$, recovering the result obtained earlier.

In summary, the hyperbolic inscribed angle theorem is related to the circular inscribed angle theorem by Wick rotation, while the parabolic inscribed angle theorem can be obtained from both by taking a limit in which the hyperbola/circle degenerates into a parabola.


\section{Acknowledgements}
The author is grateful to Maciej Dunajski and Gary Gibbons for useful discussions.


\end{document}